Mathematical Proofs 101:

How Proofs Should Be Read, Written, and Taught

Christopher M. van Bommel

University of Waterloo

**Table of Contents**

## Introduction

At some point during their schooling, students are first exposed to the idea of a mathematical proof, marking a fundamental shift in the sort of mathematics to which they are accustomed. A significant amount of research has considered mathematical proofs, the students who learn them, and the instructors that teach them, from a variety of perspectives. This paper considers this topic from four main perspectives: students' perceptions of mathematical proofs, instructors' presentations of mathematical proofs, using peer review to develop students' abilities to read proofs more critically and write proofs more convincingly, and providing students with the skills required to independently read and write proofs.

## Literature Review

### Students' Perceptions

To understand the best ways to teach mathematical proofs in an introductory course, it is important to understand how students read and process the proofs they see. According to Rav (1999), and Selden and Selden (2003), mathematicians learn new mathematics by reading and evaluating the proofs of others. This comment is supported by Weber (2004), who suggested that undergraduate mathematics students are expected to devote a significant portion of their time studying to reading proofs presented in lecture or found in textbooks. Inglis and Alcock (2012), however, noted that the skills required do not appear to be being taught successfully. As reported by Knuth (2002) and Selden and Selden (2003), both undergraduate students and precollege teachers, when asked to evaluate whether a given relatively simple argument is a valid proof, typically performed comparably to deciding by random chance. Even more discouragingly, Moore (1994) interviewed 16 mathematics and mathematics education majors, including two graduate students, and concluded "all of them said they had relied on memorizing proofs because they had not understood what a proof is nor how to write one" (p. 264).

Most of students' exposure to mathematical proofs comes during university. Uhlig (2002) observed that high school students in the United States are typically exposed to proofs in geometry but do not encounter basic proofs in other areas, for example number theory, where a basic proof would be that the sum of two odd numbers is even. Weist (2015) commented that the situation is similar in Canada. Almeida (2000) described university mathematics courses as following a "definition – theorem – proof" format, which is to say that students are taught mathematics by being presented with definitions of new concepts, theorems describing what can be derived about these concepts, and proofs giving the evidence for these derivations, and are expected to develop their own proof techniques based on the examples they see. These examples are described by Wiest as "fully realized received proofs", where the "fully realized" indicates the theorem or goal is well-stated and the proof is correct, logical, and follows a typical format with standard symbols and syntax.

Wiest (2015) detailed interviews lasting approximately ninety minutes with two students to investigate how students receive these proofs. The students chosen did not have mathematical aspirations but were hoped to be representative of undergraduate students taking mathematics classes as required in their programs. The interviews covered the participants' understanding of proofs and proving in general and in mathematics, their approaches to mathematical proofs and proving, and observed how they worked through proofs given to them. The students viewed proofs as being verifications of theorems, that is they provide evidence that the result is true and "where it came from"; as demonstrating why the result works, that is highlighting the reasoning and mechanics that goes into developing the theorem and providing a greater understanding of the ideas involved; as providing knowledge that can be used for further results; and as the result or proof can be applied. Looking at their approaches to received proofs, both students indicated using a linear method, and evaluating the proof step-by-step, but at times needed to go back when stuck or skimmed the proof first in order to get a general idea of the methods involved. The students also indicated the importance of understanding the steps of the proof for themselves and being able to recreate the reasoning, and the impact of how prior knowledge affects one's ability to understand a proof. On the other hand, neither student believed they were ever explicitly taught how to read a proof, but instead developed their approaches to reading proofs by modelling the approaches of others and learning from experience.

Inglis and Alcock (2012) also considered how students received mathematical proofs by comparing the approaches of experts and novices when reading proofs. They examined judging proof validity, students' difficulties with validating proofs, and how successful proof validations are conducted. Mejía-Ramos and Inglis (2009) argued that proof validation is distinct from reading a proof for comprehension, as the proof is already assumed to be correct. Selden and Selden (2003) argued that "construction or producing proofs is inextricably linked to the ability to validate them reliably, and a proof that could not be reliably validated would not provide much of a warrant" (p. 9). Moreover, Selden and Selden found that students tend to focus on "surface features", such as algebraic notations and computations; this finding was supported by Healy and Hoyles (2000) who found that students believed algebraic arguments, despite being "mathematically nonsensical", would receive the highest grades. Inglis and Alcock concluded that students' difficulties with validation proofs comes from overvaluing algebraic manipulations and not considering the relationships between various parts of the proof, and noted that this idea is consistent with Österholm's (2005) result that students better understood an introduction to group theory presented in words rather than both words and symbols. Weber and Mejía-Ramos (2011) observed two main strategies used to read proofs: zooming in, which is essentially the process of reading through the proof step-by-step and ensuring each step logically follows from the previous steps, filling in the missing details when necessary, and zooming out, which is essentially looking at the main steps to consider the outline of the steps of the proof.

     Inglis and Alcock (2012) looked to support or refute these ideas by tracking the eye movements of the participants in their study, both undergraduate students and academic mathematicians, as they determined the validity of given mathematical arguments, some of which had clear mistakes, and others for which the validity was debatable. For the arguments with clear mistakes, all of the mathematicians correctly identified these proofs as invalid, while the undergraduates did so about half of the time. For the other arguments, there were no significant differences observed between the two groups. The study also supported that undergraduates spend more of their time focused on formulae when reading proofs, whereas mathematicians spend more time evaluating whether the next step follows from the previous. Neither group, however, was observed to use a zooming out strategy. Inglis and Alcock noted it is important to recognize the difference in how mathematicians think they evaluate proofs compared to how they are observed to evaluate proofs, as this will aid in passing on the skills to undergraduates.

**Presenting Proofs**

Students seem to have great difficulty with mathematical proofs, despite doing well in an introduction to proofs course. Talbert (2015) noted that such a course typically marks a transition from primarily computational mathematical tasks to primarily conceptual mathematical tasks. Moreover, such courses tend to require the development of "self-regulated learning" or "independent learning" behaviours that tend not to be practised in previous mathematics courses. Pintrich (2004) formulated the model of self-regulated learning; Talbert described the model as follows:

> [S]elf-regulating learners are active participants in the learning process; have the ability to monitor and control aspects of their cognition, motivation, and behaviors related to learning; have criteria against which they can judge whether their current learning status is sufficient or whether more learning needs to take place; and that a learner's self-regulatory activities influence academic achievement. (p. 615)

These elements are typically absent in previous mathematical courses which "rely mainly on efficient computation of 'right answers'" (Talbert, 2015, p. 615). The development of independent learning as a requirement of an introduction to proofs course is supported by Lai, Weber, and Mejía-Ramos (2012), who found that proofs presented in lecture were often modified with this goal in mind, for example by leaving gaps in the proof for the students to fill. One has to be careful using this model, however, as students tend to use the proofs presented in lecture as a model for their own work (Fukawa-Connelly, 2014).

Therefore, to consider the teaching methods that are most conducive to student success, two distinct approaches to teaching are examined. Fukawa-Connelly (2014) presented a case study of an instructor's presentation of proofs in a traditional lecture format. In analyzing the proofs presented, it is important to note the distinction between "proofs that convince", that is proofs that verify a result is true, and "proofs that explain", that is proofs that provide an intuitive understanding of the correctness of the result (Hanna, 1990; Hersh, 1993). Fukawa-Connelly used Toulmin's (1969) model of argumentation to analyze the instructor's presentation of proofs in this case study and compared the level of detail in the instructor's written and spoken presentations. Toulmin's model breaks an argument into different types of statements: data – the information on which the argument is build, conclusion – the desired result, warrant – justification of the link between the data and the conclusion, and backing – supporting the validity of the warrant. The observations of Fukawa-Connelly's case study found the instructor nearly always wrote the data and the conclusion, but warrants, which are arguably the most important aspects of proofs, and backing were generally left unwritten, but sometimes stated aloud by the instructor or a student in response to the instructor's question. More precisely, warrants or backing were included written or spoken in half of the instructor's proofs, and written in only half of these. Providing the warrant or backing aloud seems to indicate the instructor saw pedagogical value in these steps, however the inconsistency in including them in proofs presented does not give students the best model for their own work.

On the other hand, Talbert (2015) investigated the benefits of an inverted classroom design in a Transition-to-Proof course for second-year mathematics majors, a course developed to teach reading and writing of mathematical proofs. Talbert notes some of the downsides to a lecture-based approach are that lecturing unintentionally creates student dependency on the lecture and that students are asked to do the more complex work of the course, that is, the actual construction of proofs, on their own time and with limited access to the instructor, which can be discouraging to those students who need assistance. With the inverted classroom model, students are first exposed to new material on their own time, through the use of readings, videos, and guided activities, and then the class time is used to allow students to work on problems with the instructor immediately available to provide assistance. Talbert used his own sections of the course as a case study, and noted most students had some difficulty adjusting to taking the course in this style, but based on student feedback, the majority had a positive view of the design and thought it to be worth the effort. He also noted that some of the challenges of implementing this structure is the large amount of time required to create materials for students to learn the material outside of the classroom, students may not be willing to adapt to the departure from the way their courses are typically taught, and that during class time the instructor must be ready to interact with students at a wide array of abilities and preparedness.

**Peer Review**

Another method of developing students' skills in writing and validating proofs is having them conduct peer reviews. Ernst, Hodge, and Schultz (2015) gave the following description of writing proofs:

> To write proofs well, students must become skeptical consumers of their own work. They must learn how to look at their own writing as if it were written in someone else's hand and ask: Does it provide the correct level of detail? Can I understand it line by line, paragraph by paragraph, and page by page? Has it correctly used and cited definitions and theorems? *Am I convinced?* (p. 121)

Zerr and Zerr (2011) commented that peer reviewing requires students to read proofs more actively compared to if they already know that the proof is correct, and would allow them to develop skills to determine if there own work is correct. Ernst, Hodge, and Schultz viewed the peer review process as providing students the opportunity to analyze the validity of arguments, examples of good proof writing to model or bad proof writing to avoid, and an audience of another student in the class rather than the instructor, which was hoped to have students focus on ensuring their proofs were understood by the reader.

Zerr and Zerr (2011) described a peer review process in which students were given a weekly "peer review problem", handing in one copy to their instructor for grading as usual, and another copy was given to a student for review. Each review was completed according to a "Peer Review Worksheet", which guided a detailed review of the proof provided. A copy was given to the author of the proof, while a second copy was submitted to the instructor for grading. Then the students were given the opportunity to correct their proofs according to the feedback they received, which required them to critically evaluate the review itself; to this end, they were required to complete an "Author's Response Sheet" which asked if they were revising their original proof, if there were suggestions given by their reviewer they were not implementing and why, and if they were making any significant changes that were not prompted by the review. The goals of this process were for the students to read proofs more actively, read proofs more critically, and write more refined proofs.

In a case study of 105 peer review cycles conducted by Zerr and Zerr (2011), 53 consisted of a correct first draft, and in only two of these were suggestions made in the review that would have led to an incorrect version. Moreover, in 19 cases, the reviewers made appropriate suggesting regarding wording or style. On the other hand, for the 52 incorrect first drafts, 38% were deemed to be correct, and only 35% had their problems completely addressed by the reviewer. Zerr and Zerr theorized that a reason for this discrepancy could be that students are more accustomed to seeing correct proofs, and had not previously been in the position to critique them. However, with all of the students completing the course earning a grade on the reviewing portion of at least 80%, Zerr and Zerr concluded that the students met the goal of becoming critical readers of proofs. This critical reading ability was also demonstrated in the students' responses to reviews and their ability to make improvements to their own work. Zerr and Zerr estimated that in 80 – 90% of the cases, the peer review process resulted in some positive learning outcome, with direct improvements of incorrect originals in 40% of those cases.

Ernst, Hodge, and Schultz (2015) also implemented a peer review exercise in two of the authors' classes, which included both within class peer review and between class peer review. Peer reviews were to identify the errors in the proof as well as comment on the overall impression the proof gave, such as systematic flaws that needed addressing or positive points to preserve, both mathematically and stylistically. The evaluations themselves were also graded by the instructor. Most students believed writing the peer reviews was more beneficial to them than receiving reviews from others, and the process made them better able to assess their own work critically. Ernst, Hodge, and Schultz commented that this believe was not surprising, since the students were going through their first experience of writing peer reviews, but noted that the quality of the reviews improved significantly during the course. The authors noted that when implementing this exercise in future courses, they would like to have them occur more frequently during the course as a way to better emphasize their importance, as well as give examples of peer reviews at the start of the course to better prepare students for writing their own.

**Proof Comprehension**

In order to write excellent proofs, it is important for students to understand the proofs to which they are exposed and to measure students' ability to comprehend proofs. Mejía-Ramos, Fuller, Weber, Rhoads, and Samkoff (2012) proposed the following seven dimensional theoretical model of proof comprehension, as presented by Hodds, Alcock, and Inglis (2014):

- Meaning of terms and statements – Understanding the meaning of symbols, terms, and definitions.
- Justification of claims – Understanding how new assertions in the proof follow from previous ones.
- Logical structure – Understanding the logical relationship between lines or components of a proof.
- Higher-level ideas – Identifying a good summary of the overarching approach of the proof.
- General method – Applying the methods within the proof to a different context.
- Application to examples – Using the ideas in the proof in terms of a specific example.
- Identifying modular structure – Understanding the main components and modules within a proof and the logical relationship between them.

Different ideas are then suggested for improving comprehension. Leron (1983) suggested the use of a "structured proof" rather than proceeding in a linear fashion, that is arranging the proofs into level with the main ideas at the top, and providing further details and justifications and subsequent levels for the steps given previously. This approach focuses on the concepts of "higher-level ideas" and "identifying modular structure", but has the drawback of separating claims from their supporting evidence. Indeed, Fuller et al. (2011) determined that students reading "structured proofs" are better able to discuss the key ideas of the proof, but had a slight drop in other aspects of proof comprehension.

Rowland (2001) suggested a generic approach, by using a specific example throughout the proof rather than providing the full argument with all the necessary algebraic notation. Such an approach is viewed as an "application to an example". This approach, however, has multiple drawbacks. Rowland found that some students did not realize that the proofs were generic and were not sufficient to prove the entire result, while others attempted to generalize the result without fully considering which relationships still applied in the more general case and which were limited to the specific example.

On the other hand, an entirely alternative approach is to provide "self-explanation training" to the students, the idea that students create their own explanations of the proofs they are required to learn; the term is due to Chi et al. (1989). A study was conducted by Hodds, Alcock, and Inglis (2014) in which students were taught the basic principles of self-explanation training: "identifying key ideas in each line of a proof, and explaining each line in terms of previous ideas presented in the proof or in terms of previous knowledge" (p. 22), and these students were compared to a control group in assessing their success at comprehending a given proof. They found that the student that were given self-explanation training were able to give better quality explanations of the proof and performed significantly better on the comprehension test. In a longer term study, they obtained results that self-explanation training also has lasting effects. Interestingly, there was not a significant increase in time spent reading the proofs for the students given self-explanation training compared to those who had not, which seems to indicate the students are making better use of their time spent reading the proof to obtain a fuller understanding.

## Discussion

As indicated by Uhlig (2002) and Wiest (2015), students coming into university are, for the most part, seeing proofs for the first time, whether it be in a course specifically designed to introduce proofs, or a course for which proofs are an essential part of the material, for example, algebra, combinatorics, or number theory. While some students are exposed to geometric proofs during high school, these proofs tend to follow a very different structure. Such proofs tend to take a two column approach, where statements in the first column are the facts and claims that make up the proof, and information in the second column gives the rationale for such a statement, whether it be given as part of the problem, or the statement follows from a particular result or theorem that was given previously. On the other hand, proofs taught in university tend to be more expository, structured in sentences and paragraphs, with a mixture of mathematical notation and algebraic expressions. It is a significant departure from the "show your work" mentality of high school and introductory calculus courses, as suggested by Zerr and Zerr (2001), where listing a series of mathematical equations is typically sufficient to earn full marks, although this request does go some way to demonstrating that the entire process is important rather than just obtaining the correct answer.

It therefore seems that the majority of students are unprepared to make this jump. Most courses do not seem to illustrate what makes a good proof or how students can read, understand, and write proofs on their own, but instructors instead appear to expect students to learn from their examples. As indicated by the case study of Fukawa-Connelly (2014), however, instructors do not necessarily model best practices in their classes. Moreover, it was indicated from the study of Inglis and Alcock (2012) that mathematicians cannot even always accurately describe their own process for reading and understanding proofs, as the observations from their eye-movement study tended to contradicts the mathematicians descriptions of their strategy. Hence, despite their best intentions, instructors may not be adequately preparing students for the process of reading and writing proofs, because they do not themselves realize they are not following their own advice.

It may also be the case that instructors feel pressured to ensure they cover an adequate amount of material in their classes and so find it necessary to skip steps of proofs when writing on the board, even if they supply the necessary step orally, or skip the step altogether it a similar sort of proof was done previously. Hence, instructors should take care to write about proofs fully, especially at the beginning of introductory courses involving proofs, so that their students do not get the wrong idea of what is sufficient and what can be left out. Better still, it would be beneficial for instructors of these courses to more openly discuss, particularly early on, how each element of a proof contributes to the presentation and whether it would be as convincing without.

Having students write peer reviews for each other's work seems to immediately follow on from this idea. What better way to demonstrate to students the flaws in their own proofs by having them figure out if they are convinced by someone else's argument? Peer reviewing provides examples that are rarely typical in textbooks or the instructor's lecture notes because they have the potential to be incorrect. In some ways, identifying the errors present demonstrates an even better understanding of the result than proving the result in the first place, as it is something that cannot simply be achieved by modelling a similar proof.

Using peer reviews in courses also introduces students to a reality that exists in academia. Any result published in a mathematical journal, or any other scholarly journal, goes through a peer review process, where the reviewer must be convinced of the correctness of the result, and may also make suggestions regarding style and presentation. Thinking about presenting mathematics in this way that a complete stranger would be able to follow it may be an important factor for undergraduate students, who are used to consider writing for an audience that knows the material they have been taught in the course; most often thinking of the instructor. Writing proofs with others in mind should only improve a student's ability to produce elegant proofs.

So in some ways, given the way these courses are currently structured, the instructors are doing their students a great disservice. They are expected to read and write proofs on their own without ever being given the skills to do so. Instead, it seems to be hoped they will learn by exposure and repetition. Unfortunately, as evidenced by Knuth (2002) and Selden and Selden (2003), students may as well guess at whether a given argument is correct, and discouragingly, students are reported to resort to memory for producing proofs (Moore, 1994).This circumstance may be partially a result of the testing students are subjected to; because they know they will have to produce either the same or similar proofs to those they have seen in class or on assignments, it is simpler for them to memorize the content and structure of the proofs, rather than understanding or appreciating them and hence being able to extend their knowledge to new ideas.

Hence, this idea of "self-explanation training" presented by Hodds, Alcock, and Inglis (2014), is very encouraging. It demonstrates a simple strategy that could be added to any introductory mathematics course that would give students the skills they need to read and write proofs. It seems it is important to demonstrate to students that reading proofs is not a passive activity, but that while reading, every statement should be questioned to ensure understanding of why each step is valid, and how the argument works as a whole. Ideally, students will also question their own proofs in a similar manner, so they can decide for themselves whether it is sufficient and convincing, or if more explanation is necessary.

## Conclusion

Students who are first exposed to mathematical proofs often do not seem to have this exposure accompanied by a development of the skills required to read proofs with a goal of comprehension and write proofs with a goal of creating a self-contained convincing argument. Instructors should recognize students may need more guidance when they are first exposed to proofs than they can obtain by simply following the examples shown in class or the textbook and doing their best to replicate those examples. In fact, simply replicating the examples presented does not demonstrate a full understanding of the proofs they are required to produce in assignments or exams. By understanding students' perceptions to proofs, taking care in what behaviour is modelled in the classroom or even adapting the classroom structure to better accommodate students' needs, and providing opportunities for students to think about proofs critically, instructors should be able to more effectively teach students how to read, write, understand, and evaluate the proofs they will encounter throughout their courses.